\documentclass[letterpaper, 10pt, conference]{ieeeconf}
\usepackage[letterpaper, top=57pt, left=48pt, right=54pt, bottom=45pt]{geometry}

\IEEEoverridecommandlockouts

\usepackage{graphicx}
\usepackage{grffile}
\usepackage{amsmath}
\usepackage{amssymb}
\usepackage{tikz}
\usepackage[font=small]{caption}
\usepackage{subcaption}
\usepackage{pgfplots}
\usepackage{dblfloatfix}
\usepackage{float}
\usepackage[linesnumbered,noend]{algorithm2e}
\SetAlCapNameFnt{\small}
\SetAlCapFnt{\small}

\usepackage{cuted}
\usetikzlibrary{arrows.meta, calc, decorations.pathreplacing, decorations.markings}
\usepgfplotslibrary{polar}

\tikzset{
    boat/.pic={
        \draw[draw=black] (-0.3, 0) -- (0.3, 0);
        \draw[draw=black] (-0.3, 0) parabola[bend at end] (0, 1.0);
        \draw[draw=black] (0.3, 0) parabola[bend at end] (0, 1.0);
    }
}

\definecolor{darkgreen}{rgb}{0.25,0.75,0.25}
\definecolor{softred}{rgb}{1.0,0.3,0.3}
\definecolor{mygreen}{HTML}{5BAD2B}
\definecolor{myblue}{HTML}{2B5BAD}
\definecolor{mypurple}{HTML}{7D2BAD}
\definecolor{mypink}{HTML}{AD2B5B}
\definecolor{myorange}{HTML}{F39831}
\colorlet{windgreen}{green!60!black}
\AtBeginDocument{\colorlet{defaultcolor}{.}}

\usepackage{url}
\usepackage[pdfstartview=FitV, urlcolor=blue]{hyperref}

\newcommand{\bx}{\boldsymbol{x}}
\newcommand{\by}{\boldsymbol{y}}

\newcommand{\br}{\boldsymbol{r}}

\newcommand{\bp}{\boldsymbol{p}}

\newcommand{\x}{\bx}
\newcommand{\y}{\by}

\DeclareMathOperator{\atantwo}{atan2}

\mathchardef\mhyphen="2D

\title{\LARGE \bf
Stochastic Optimal Control of a Sailboat
}

\author{Cole Miles$^{1,2}$, Alexander Vladimirsky$^{3}$%
\thanks{*Supported in part by the NSF DMS (awards 1738010 and 2111522), and the DOE CSGF (award DE-SC0020347).}%
\thanks{$^{1}$Department of Physics, Cornell University}%
\thanks{$^{2}$cmm572@cornell.edu}%
\thanks{$^{3}$Department of Mathematics, Cornell University}%
}

\begin{document}

\maketitle
\thispagestyle{empty}
\pagestyle{empty}

%%%%%%%%%%%%%%%%%%%%%%%%%%%%%%%%%%%%%%%%%%%%%%%%%%%%%%%%%%%%%%%%%%%%%%%%%%%%%%%%
\begin{abstract}

% !TEX root = root.tex
In match race sailing, competitors must steer their boats upwind in the presence of unpredictably evolving weather. 
Combined with the tacking motion necessary to make upwind progress, this makes it natural to model their path-planning as a hybrid stochastic optimal control problem.
Dynamic programming provides the tools for solving these, but the computational cost can be significant.
We greatly accelerate a semi-Lagrangian iterative approach 
of Ferretti and Festa \cite{ferretti_optimal_2019} 
by reducing the state space dimension and designing an adaptive timestep discretization that is very nearly causal. We also provide a more accurate tack-switching operator by integrating over potential
wind states after the switch.
The method is illustrated through a series of simulations with varying stochastic wind conditions.

\end{abstract}

%%%%%%%%%%%%%%%%%%%%%%%%%%%%%%%%%%%%%%%%%%%%%%%%%%%%%%%%%%%%%%%%%%%%%%%%%%%%%%%%
\section{INTRODUCTION}

% !TEX root = root.tex
In recent years, advances in numerical and mathematical methods have greatly impacted the sport of sailing. Yacht path planning algorithms are now heavily employed both as a part of match simulation software, and as a tool to provide live guidance on optimal decision-making during a race. Here, we focus on the problem of finding the optimal feedback control to perform live path planning in the presence of a stochastically-evolving weather system. 

Weather conditions greatly impact the performance capabilities of a sailboat in a race. While it is possible to produce meteorological forecasts, there is an inevitable stochastic component to how weather conditions evolve during the course of a race which is impossible to predict exactly. It is critical that algorithms used for decision making correctly factor this randomness into their control policy. In \cite{philpott_simulation_2004}, this randomness is modeled as a discrete-time Markov chain process which updates weather conditions at a fixed tick rate. In \cite{vinckenbosch_stochastic_2012}  this idea is extended to a continuous-time Markov chain, formulating the system as a free boundary problem whose solution is a switching curve giving the optimal control policy. Here, we will take the approach of \cite{ferretti_optimal_2019} and describe optimality through dynamic programming, yielding a system of two coupled quasi-variational inequalities of Hamilton-Jacobi-Bellman (HJB) type with degenerate ellipticity.

For full physical accuracy, a model would have to
include both the position and velocity of the sailboat in the state space, and model the forces on the sailboat in different state configurations, leading to at best $5$ coupled first-order equations in $5$ state variables. An approach along these lines can be found in \cite{philpott_simulation_2004}. But in long course matches, inertial effects associated with accelerating up to the equilibrium speed are (in general) secondary, allowing one to reduce the system to $3$ coupled first-order equations by assuming that the velocity vector can be changed instantaneously.

However, this assumption cannot be made in modeling
the process of \textit{tacking}, which is necessary when attempting to travel upwind. To make progress, sailors must travel at angles to the wind, repeatedly switching directions to steer in a zig-zag pattern on ether side of the wind direction. Each time the bow is swung to the other side of the wind (a single \textit{tack}), the boat is significantly slowed down for the period of time in which the wind pushes against the boat. In previous works (\cite{ferretti_optimal_2019},\cite{cacace_stochastic_2019}), this is incorporated by modeling the system as a \textit{hybrid control} problem, with a continuous control corresponding to the angle of the sails and a discrete control corresponding to which tack we should be on. The lost inertia has been modeled as a fixed switching delay, reflecting the time lost due to the reduced speed from performing the tack.

In this work, we extend and improve upon the methods of Ferretti and Festa \cite{ferretti_optimal_2019}, who solve the problem using a semi-Lagrangian discretization. In Sec.~\ref{s:problem_setup}, we give the original problem statement and recast it to a reduced coordinate system in 2D, which alone greatly reduces the time needed to solve the problem numerically. In Sec.~\ref{s:theory}, we set up the optimal control problem following the structure of Ref.~\cite{ferretti_optimal_2019}, and provide an improvement to their \textit{switching operator} to incorporate the wind evolution during the process of a tack. In Sec.~\ref{s:numerics}, we describe our semi-Lagrangian numerical methods, including a discretization in which the system is nearly-causal and can be solved in a handful of iterations. We present the results of numerical simulations in Sec.~\ref{s:experiments}, producing optimal switching maps and highlighting the effect of our improved switching operator. We conclude by listing directions for future work in Sec.~\ref{s:conclusions}.

\section{System Dynamics}
\label{s:problem_setup}
% !TEX root = root.tex
\begin{figure*}[]
    \centering
    \begin{subfigure}{0.38\textwidth}
        \centering
        \begin{tikzpicture}
            \node (polarplot) at (0, 0) {
                \resizebox{0.82\textwidth}{!}{
                    \begin{tikzpicture}[]
\begin{polaraxis}[height = {127.0mm}, ylabel = {}, xmin = {90}, xmax = {450}, ymax = {0.05500000000000001}, xlabel = {}, unbounded coords=jump,scaled x ticks = false,xlabel style = {font = {\fontsize{11 pt}{14.3 pt}\selectfont}, color = {rgb,1:red,0.00000000;green,0.00000000;blue,0.00000000}, draw opacity = 1.0, rotate = 0.0},xmajorgrids = true,xtick = {90.0,135.0,180.0,225.0,270.0,315.0,360.0,405.0},xticklabels = {$90$,$135$,$180$,$225$,$270$,$315$,$0$,$45$},xtick align = inside,xticklabel style = {font = {\fontsize{8 pt}{10.4 pt}\selectfont}, color = {rgb,1:red,0.00000000;green,0.00000000;blue,0.00000000}, draw opacity = 1.0, rotate = 0.0},x grid style = {color = {rgb,1:red,0.00000000;green,0.00000000;blue,0.00000000},
draw opacity = 0.1,
line width = 0.5,
solid},axis x line* = left,x axis line style = {color = {rgb,1:red,0.00000000;green,0.00000000;blue,0.00000000},
draw opacity = 1.0,
line width = 1,
solid},scaled y ticks = false,ylabel style = {font = {\fontsize{11 pt}{14.3 pt}\selectfont}, color = {rgb,1:red,0.00000000;green,0.00000000;blue,0.00000000}, draw opacity = 1.0, rotate = 0.0},ymajorgrids = true,ytick = {0.0,0.01,0.02,0.03,0.04,0.05},yticklabels = {$0.00$,$0.01$,$0.02$,$0.03$,$0.04$,$0.05$},ytick align = inside,yticklabel style = {font = {\fontsize{8 pt}{10.4 pt}\selectfont}, color = {rgb,1:red,0.00000000;green,0.00000000;blue,0.00000000}, draw opacity = 1.0, rotate = 0.0},y grid style = {color = {rgb,1:red,0.00000000;green,0.00000000;blue,0.00000000},
draw opacity = 0.1,
line width = 0.5,
solid},axis y line* = left,y axis line style = {color = {rgb,1:red,0.00000000;green,0.00000000;blue,0.00000000},
draw opacity = 1.0,
line width = 1,
solid},    xshift = 0.0mm,
    yshift = 0.0mm,
    axis background/.style={fill={rgb,1:red,1.00000000;green,1.00000000;blue,1.00000000}}
,legend style = {color = {rgb,1:red,0.00000000;green,0.00000000;blue,0.00000000},
draw opacity = 1.0,
line width = 1,
solid,fill = {rgb,1:red,1.00000000;green,1.00000000;blue,1.00000000},fill opacity = 1.0,text opacity = 1.0,font = {\fontsize{8 pt}{10.4 pt}\selectfont}},colorbar style={title=}, ymin = {0}, width = {127.0mm}]\addplot+ [color = {rgb,1:red,0.00000000;green,0.60560316;blue,0.97868012},
draw opacity = 1.0,
line width = 1,
solid,mark = *,
mark size = 2.0,
mark options = {
            color = {rgb,1:red,0.00000000;green,0.00000000;blue,0.00000000}, draw opacity = 1.0,
            fill = {rgb,1:red,0.00000000;green,0.60560316;blue,0.97868012}, fill opacity = 1.0,
            line width = 1,
            rotate = 0,
            solid
        },forget plot]coordinates {
(0.0, 0.0)
(5.0, 0.004615384615384616)
(10.0, 0.009230769230769232)
(14.999999999999998, 0.013846153846153847)
(20.0, 0.016153846153846154)
(25.0, 0.019230769230769232)
(32.0, 0.031538461538461536)
(36.0, 0.03538461538461538)
(40.0, 0.038461538461538464)
(45.0, 0.041538461538461545)
(52.0, 0.04461538461538462)
(59.99999999999999, 0.047692307692307694)
(70.0, 0.04923076923076924)
(80.0, 0.05)
(90.0, 0.05)
(100.0, 0.05)
(110.0, 0.04923076923076924)
(119.99999999999999, 0.047692307692307694)
(130.0, 0.04384615384615385)
(140.0, 0.04)
(150.0, 0.036153846153846154)
(160.0, 0.03230769230769231)
(170.0, 0.030000000000000002)
(180.0, 0.027692307692307693)
};
\addplot+ [color = {rgb,1:red,1.00000000;green,0.00000000;blue,0.00000000},
draw opacity = 1.0,
line width = 1,
solid,mark = *,
mark size = 2.0,
mark options = {
            color = {rgb,1:red,0.00000000;green,0.00000000;blue,0.00000000}, draw opacity = 1.0,
            fill = {rgb,1:red,1.00000000;green,0.00000000;blue,0.00000000}, fill opacity = 1.0,
            line width = 1,
            rotate = 0,
            solid
        },forget plot]coordinates {
(-0.0, 0.0)
(-5.0, 0.004615384615384616)
(-10.0, 0.009230769230769232)
(-14.999999999999998, 0.013846153846153847)
(-20.0, 0.016153846153846154)
(-25.0, 0.019230769230769232)
(-32.0, 0.031538461538461536)
(-36.0, 0.03538461538461538)
(-40.0, 0.038461538461538464)
(-45.0, 0.041538461538461545)
(-52.0, 0.04461538461538462)
(-59.99999999999999, 0.047692307692307694)
(-70.0, 0.04923076923076924)
(-80.0, 0.05)
(-90.0, 0.05)
(-100.0, 0.05)
(-110.0, 0.04923076923076924)
(-119.99999999999999, 0.047692307692307694)
(-130.0, 0.04384615384615385)
(-140.0, 0.04)
(-150.0, 0.036153846153846154)
(-160.0, 0.03230769230769231)
(-170.0, 0.030000000000000002)
(-180.0, 0.027692307692307693)
};
\end{polaraxis}

\end{tikzpicture}    
                }
            };        
            \draw[-{Stealth}, very thick, color=windgreen] (4, 0) -- node[above]{\color{windgreen}Wind} (0, 0);
            \node at (0.1, 1.6) {\color{blue!80!green}Tack $q = 1$};
            \node at (0.1, -1.6) {\color{red}Tack $q = 2$};
        \end{tikzpicture}

        \caption{The polar speed plot $s(u)$ used in this work. Note the symmetry about $u = 0^\circ$, non-convexities near $0^\circ$ and $180^\circ$, and that the speed vanishes at $u = 0^\circ$. These properties are common to all polar plots.}
        \label{fig_polarplot}
    \end{subfigure}
    \begin{subfigure}{0.61\textwidth}
        \centering
        \begin{tikzpicture}[thick,scale=0.85, every node/.style={scale=0.85}]
            % Tack label
            \node (tacklabel) at (-1, 5) {Tack $q = 1$};

            % Target circle
            \draw[thick, draw=myblue, fill=myblue!20] (0, 4) node(target){$\mathcal{T}$} circle [radius = 0.5];
            % Boat
            \draw (0, 0) node(boat){} pic[rotate=75] {boat};

            % Tack arc
            \draw[draw=black, fill=blue, opacity=0.5, shift=(140:0.5)] (0, 0) arc (140:320:0.5);

            % Radial distance
            \path[draw, decorate,decoration={brace,mirror,amplitude=10pt}, xshift=14pt] (0, 0) -- node[right=12pt]{\color{black}$r$} (0, 4);

            % Wind
            \draw[-{Stealth}, thick, draw=windgreen, rotate around={50:(0,0)}] (0, 4) node[above right] {\color{windgreen}Wind} -- (0,-1);

            % Vertical guide line
            \draw[dashed] (0, 4) to (0, 0);
            % Boat speed vector
            \draw[-{Stealth}, thick, draw=mypurple, rotate=75] (0, 0) -- (0, 3.5) node[below left] {\color{mypurple}$\vec{s}$} node (speedend) {};

            % Angle arcs
            \draw[draw=mygreen, rotate=0, ->] (0, 2.0) arc (90:140:2.0) node[above=5pt] {\color{mygreen}$\theta$};
            \draw[draw=mypurple, rotate=50, ->] (0, 1.3) arc (90:115:1.3) node[above=5pt] {\color{mypurple}$u$};

            % Speed components
            \draw[-{Stealth}, draw=mypurple!60] let \p1 = (speedend) in (0, 0) -> (0, \y1) node [below right] {\color{mypurple!60}$dr$};
            \draw[-{Stealth}, draw=mypurple!60] let \p1 = (speedend) in (0, \y1) -> (\x1, \y1) node [above right=5pt] {\color{mypurple!60}$d\theta$};

        \end{tikzpicture}
        \begin{tikzpicture}[thick,scale=0.85, every node/.style={scale=0.85}]
            % Tack label
            \node (tacklabel) at (-1, 5) {Tack $q = 2$};

            % Target circle
            \draw[thick, draw=myblue, fill=myblue!20] (0, 4) node(target){$\mathcal{T}$} circle [radius = 0.5];
            % Boat
            \draw (0, 0) node(boat){} pic[rotate=25] {boat};

            % Tack arc
            \draw[draw=black, fill=red, opacity=0.5, shift=(140:0.5)] (0, 0) arc (140:-40:0.5);

            % Radial distance
            \path[draw, decorate,decoration={brace,mirror,amplitude=10pt}, xshift=14pt] (0, 0) -- node[right=12pt]{\color{black}$r$} (0, 4);

            % Wind
            \draw[-{Stealth}, thick, draw=windgreen, rotate around={50:(0,0)}] (0, 4) node[above right] {\color{windgreen}Wind} -- (0,-1);

            % Vertical guide line
            \draw[dashed] (0, 4) to (0, 0);
            % Boat speed vector
            \draw[-{Stealth}, thick, draw=mypurple, rotate=25] (0, 0) -- (0, 3.5) node[below left] {\color{mypurple}$\vec{s}$} node (speedend) {};

            % Angle arcs
            \draw[draw=mygreen, rotate=0, ->] (0, 2.0) arc (90:140:2.0) node[above=5pt] {\color{mygreen}$\theta$};
            \draw[draw=mypurple, rotate=25, <-] (0, 1.3) arc (90:115:1.3) node[above=5pt] {\color{mypurple}$u$};

            % Speed components
            \draw[-{Stealth}, draw=mypurple!60] let \p1 = (speedend) in (0, 0) -> (0, \y1) node [below left] {\color{mypurple!60}$dr$};
            \draw[-{Stealth}, draw=mypurple!60] let \p1 = (speedend) in (0, \y1) -> (\x1, \y1) node [above right=5pt] {\color{mypurple!60}$d\theta$};
        \end{tikzpicture}
        \caption{A diagram of the system setup in reduced coordinates.}
        \label{fig_reduced_setup}
    \end{subfigure}
    \caption{Boat dynamics relative to the wind and relative to a target.}
\end{figure*}
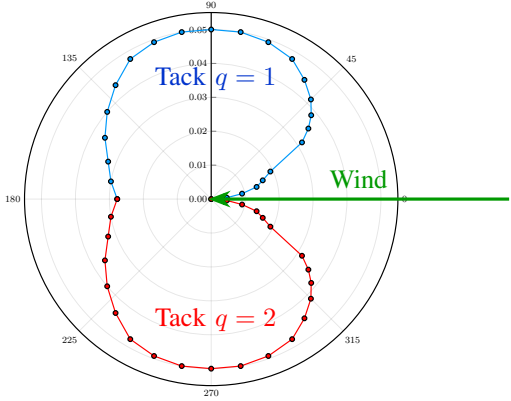
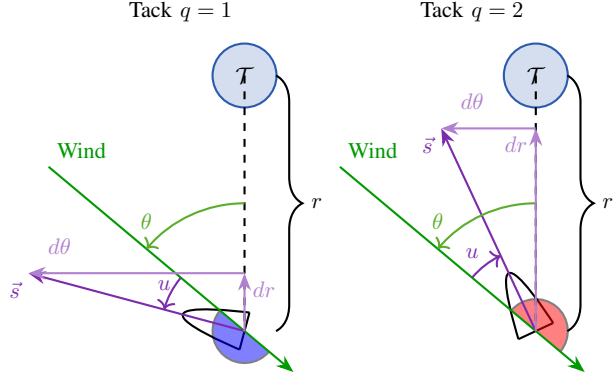

We consider a problem introduced in \cite{ferretti_optimal_2019}:
time-optimal sailboat navigation to a closed target set $\mathcal{T}$ when the wind direction evolves stochastically. For simplicity, we focus on their simplified model where the wind speed is fixed, and the angle $\phi$ of the upwind direction (measured counterclockwise from the $y$-axis) evolves under a drift-diffusion process:
\begin{equation}
    d\phi = a dt + \sigma dW
\end{equation}
where $a$ is the drift rate, $\sigma$ is the strength of the Brownian motion, and $dW$ is the differential of a Wiener process.

We control the angle of the boat's velocity relative to the upwind direction.
But at any given time we are also restricted to an interval of angles determined by the sailboat's current \textit{tack}, or configuration of sails. We represent the tack as a controlled discrete variable $q \in \{1, 2\}$; $q = 1$ represents the starboard tack, which we model as having access to the angles $[0, \pi]$, while the port tack $q = 2$ has access to $[-\pi, 0]$. For notational convenience, we will use $(-1)^q u$ to encode the velocity angle relative to the upwind direction,
with the steering control variable $u \in [0, \pi].$ So, the steering angle $u$ is measured counterclockwise for $q=1$ and clockwise for $q=2,$
with $u=0$ corresponding to motion directly against the wind. Due to our inertia-less approximation, whenever $u$ changes, the boat is assumed to instantly start travelling in the new direction with the angle-dependent speed $s(u)$.  The speed profile (termed the \textit{polar plot}) capturing this dependence is determined by the exact geometry of a specific boat. A typical example \cite{noauthor_orc_nodate} is provided in Fig.~\ref{fig_polarplot}.

Assuming that the boat's location is $(x(t),y(t)),$ the system dynamics in absolute coordinates is then
\begin{align}
    dx &= -s(u) \sin(\phi - (-1)^qu) dt \nonumber \\
    dy &= s(u) \cos(\phi - (-1)^qu) dt \\
    d\phi &= adt + \sigma dW \nonumber
\end{align}

We now simplify by focusing on radially symmetric problems:
if we assume that the target $\mathcal{T}$ is circular, we no longer care about an absolute angle of the boat's position relative to the target. Shifting to a coordinate system centered at $\mathcal{T},$ we can
reduce the $(x, y, \phi)$ state space to just two coordinates $(r, \theta)$, where $r = \sqrt{x^2 + y^2}$ is the radial distance of the boat from the center of $\mathcal{T}$ and $\theta  = \phi + \atantwo(-x, -y)$ is the angle of the wind measured counterclockwise relative to the straight line between the boat and the target; see Fig.\ \ref{fig_reduced_setup}. In this reduced representation, the dynamics are governed by
\begin{align}
    dr &= -s(u)\cos(\theta - (-1)^q u) dt\nonumber \\
    d\theta &= \left[\frac{s(u)}{r} \sin(\theta - (-1)^q u) + a\right]dt + \sigma dW
    \label{eq_reducedcoords}
\end{align}
which describe the projections of the velocity vector onto the radial and angular components, combined with the drift-diffusion evolution of the wind. We refer to the coefficients on $dt$ in each equation as $\dot{r}$ and $\dot{\theta}$, respectively. New from the original system is the $\dot{\theta} \propto \frac{1}{r}$ relationship, which will force us to utilize an adaptive timestepping scheme that chooses smaller timesteps near the target to counteract the rapidly increasing angular speed. We note that $\dot{\theta}$ still remains bounded since $\mathcal{T}$ has positive radius.

\section{Stochastic Optimal Control}
\label{s:theory}
% !TEX root = root.tex
Assuming we start from a position $\br = (r, \theta)$ and tack $q$, the total time-to-target $T_\mu(\br, q) = \min \{t \mid \br'(t) \in \mathcal{T}, \br'(0) = \br, q'(0) = q  \}$ is a random variable whose expectation we are trying to minimize\footnote{
Strictly speaking,
the objective in \cite{ferretti_optimal_2019} is different since they minimize $\mathbb{E} \left[ \int_0^T e^{-\lambda t} dt \right] = \mathbb{E} \left[ \left(1 - \exp(- \lambda T) \right) / \lambda \right],$ where $\lambda > 0$ is a time-discounting coefficient
chosen to improve the convergence of value iterations.  Our objective is obtained by taking the limit with $\lambda \to 0$ and the numerical tests indicate that value iterations based on our discretization converge even in this ``undiscounted'' case.
} by choosing a feedback steering/tack-switching policy $\mu$.
This $T_\mu(\br, q)$ can be split into two pieces: the time spent continuously controlling the sails plus the time spent performing tack-switches. Each tack-switch incurs a \textit{switching delay} $C$ due to the speed loss from the wind pushing directly against the boat's intended direction during the maneuver\footnote{For simplicity of modeling, we will treat this delay as constant, though state-dependent $C$ would not present additional computational challenges.
}. We model the boat's propelled speed as falling to zero (i.e., $s(u) = 0$) during a tack-switch.

We thus define the \textit{value function} to encode the minimal expected time to $\mathcal{T}$ form each starting configuration:
\begin{equation}
    v(\br, q) := \inf_\mu \mathbb{E} \left[ T_\mu(\br, q) \right].
    \label{eq_valuefunc}
\end{equation}
By Bellman's Optimality Principle, this $v$ must satisfy 
\begin{equation}
    v(\br, q) = \min\left(\tau + \inf_{\mu}\Sigma_{\mu, \tau} v(\br, q), \hspace{0.5em} C + \mathcal{N}v(\br, q)\right) + o(\tau),
    \label{eq_bellmanopt}
\end{equation}
where the first argument to the $\min$ is to be interpreted as the best expected time-to-$\mathcal{T}$ \textit{if we stay on the current tack} at least for  a small time $\tau$, while the second is the best expected time $\textit{if we switch tacks immediately}$. More formally, we define the operators $\Sigma_{\mu, \tau}$ and $\mathcal{N}$ as:
\begin{align}
    \Sigma_{\mu, \tau} v(\br, q) &= \mathbb{E}_{\br' ; \mu, \tau}[v(\br', q) | \br], \label{eq_control_op} \\
    \mathcal{N} v(\br, q) &= \mathbb{E}_{\br' ; \O, C}[v(\br', 3 - q) | \br], \label{eq_switch_op}
\end{align}
where we employ the compact notation $\mathbb{E}_{\br'; \mu,\tau}[\cdot | \br]$ to indicate the expectation of a function over evolved positions $\br'$, starting from position $\br$ and using control policy $\mu$ for a time $\tau$. In \eqref{eq_switch_op}, the symbol $\O$ in place of the control policy indicates that during the tacking time $C$, no steering control is taking place as $s(u) = 0$ for the duration.
Importantly, the version of operator \eqref{eq_switch_op} used in $\cite{ferretti_optimal_2019}$ does not treat $\theta$ as evolving during the tack, instead setting $\mathcal{N} v(\br, q) = v(\br, 3 - q)$. As we show in section \ref{s:experiments}, correctly incorporating these switch-duration dynamics is critical in high-drift problems.

A standard argument based on an Ito-Taylor expansion of \eqref{eq_bellmanopt} yields the quasi-variational inequality that must be satisfied by a smooth value function: 
\begin{equation}
    \max\left(v - \mathcal{N}v - C,  \quad H(\br, q, \nabla_{\br}^{}v) - \frac{\sigma^2}{2} v_{\theta \theta}^{} \right) = 0 \label{eq_hjb}
\end{equation}
with the Hamiltonian
\begin{equation}
    H(\br, q, \bp) = \max_u \left(-\dot{r}(\br, q, u) p_1 - \dot{\theta}(\br, q, u) p_2\right) - 1.
\end{equation}
If $v$ is not smooth, one can still interpret it as a weak solution (in particular, the unique {\em viscosity solution} \cite{BENSOUSSAN2000261}) of \eqref{eq_hjb}.

As $\theta$ is the only randomly evolving variable and obeys drift-diffusive dynamics, both expectations in \eqref{eq_control_op} and \eqref{eq_switch_op} take a similar form as averaging $v$ over Gaussian-distributed final $\theta'$. For small $\tau$, an approximation to \eqref{eq_control_op} taking $u=\mu(\br, q)$ and $\dot{r}, \dot{\theta}$ constant for a time $\tau$ yields:
\begin{multline}
    \Sigma_{u, \tau} v(r, \theta, q) = \frac{1}{\sqrt{2\pi \tau}\sigma}\int_{-\infty}^{\infty}d\theta' \\ \exp\left\{-\frac{(\theta' - \theta - \dot{\theta}\tau)^2}{2\tau \sigma^2}\right\}v(r + \dot{r}\tau, \theta', q), 
    \label{eq_control_op_integral}
\end{multline}
where $v$ is interpreted as $2\pi$-periodic in $\theta'$, and $\dot{r}, \dot{\theta}$ both implicitly depend on $(r, \theta, q, u)$.
Eq.~\eqref{eq_switch_op} is the same formula under the substitutions $q \to 3 - q, \tau \to C, \dot{r} = 0, \dot{\theta} = a$.

\section{Numerical Implementation}
\label{s:numerics}

\subsection{Semi-Lagrangian Discretization}
\label{ss:semi_lagrangian}
% !TEX root = root.tex

For a fixed small $\tau$, an approximation of the value function  $v$ can be obtained as a fixed point of the mapping defined by ignoring the $o(\tau)$ term in \eqref{eq_bellmanopt}. It can be shown that if the schemes approximating both arguments in the $\min$ are monotone, stable, and consistent with \eqref{eq_hjb}, then iterating this mapping is guaranteed to converge \cite{barles_convergence_1991}. As in \cite{ferretti_optimal_2019}, to build a monotone scheme we discretize the state space on a rectangular grid and perform $\textit{value iterations}$ of \eqref{eq_bellmanopt} starting from an initial guess $v^0$ to produce a sequence of approximate value functions $v^n$ converging to the solution of the discretized version of \eqref{eq_bellmanopt}. Using monotone schemes, choosing $v^0$ such that the first iteration produces $v^1 \le v^0$ everywhere is sufficient to ensure convergence \cite{Ferretti2015JOptimTheoryAppl}.

To approximate the continuous-control operator $\Sigma_{\mu, \tau},$ we change the variable of integration in \eqref{eq_control_op_integral} to $\xi = (\theta' - \theta - \dot{\theta}\tau)/\sqrt{2\tau\sigma^2}$ and use the Gauss-Hermite quadrature \cite{abramowitz_handbook_2013}:
\begin{multline}
    \Sigma_{u, \tau} v^n(\br, q) \\ = \frac{1}{\sqrt{\pi}} \int_{-\infty}^{\infty} d\xi e^{-\xi^2}v^n(r + \dot{r}\tau, \sqrt{2 \tau}\sigma \xi + \theta + \dot{\theta}\tau, q) \\ 
    \approx \frac{1}{\sqrt{\pi}}\sum_{i=1}^M w_i v^n(r + \dot{r}\tau, \, \sqrt{2\tau}\sigma \xi_i + \theta + \dot{\theta}\tau, \, q), \label{eq_gausshermite_cont}
\end{multline}
where $M$ is the number of quadrature points used, $\xi_i$ are the roots of the $M$\textsuperscript{th} Hermite polynomial $H_M(\xi)$, and $w_i$ are 
\begin{equation}
    w_i = \frac{2^{M-1}M!\sqrt{\pi}}{M^2[H_{M-1}(\xi_i)]^2}.
\end{equation}

The sum in \eqref{eq_gausshermite_cont} requires evaluating $v^n$ off-grid, which we accomplish through interpolation, as is common for \textit{semi-Lagrangian schemes}~\cite{Falcone_book}.  The latter are guaranteed to be monotone as long as the interpolation used is zeroth or first order \cite{ferretti_technique_2010}. In our implementation, we define the numerical operator $S_{u, \tau}$ by the sum of \eqref{eq_gausshermite_cont} evaluated using bilinear interpolation between the 4 closest gridpoints. 
We use $M = 2$ quadrature points, which reduces $S_{u, \tau}$ to the form in \cite{ferretti_optimal_2019}, averaging between the two points 
$\br_{\pm} = \left( r + \dot{r} \tau, \, \theta + \dot{\theta} \tau \pm \sigma \sqrt{\tau} \right).$

Due to $\dot{\theta}$ rapidly increasing as $\frac{1}{r}$ near the target, we adaptively choose our time discretization as:
\begin{equation}
    \tau(\br, q, u) \propto \min\left(\frac{\Delta r}{\dot{r}(\br, q, u)}, \frac{\Delta \theta}{\dot{\theta}(\br, q, u)}\right). \label{eq_adaptivedt}
\end{equation}
We choose the constant of proportionality in \eqref{eq_adaptivedt} to be $1.5$, to step into the middle of the first grid cell along the direction of motion. We then denote the scheme $S$ as the optimization of this operator over control angles $u$:
\begin{equation}
    S v^n(\br, q) = \min_u\Big( \tau(\br, q, u) + S_{u,\tau(\br, u, q)} v^n(\br, q) \Big).
\end{equation}
The optimal control policy solution is then precisely the $\arg\min$ of (14). In practice, since our polar plot datasets contain a discrete set of $s(u)$ values, we perform this minimization by a simple grid search over the tabulated angles.
The accuracy of this procedure can be also improved by expanding the available set of control angles through interpolation.

The same reparameterization and Gauss-Hermite quadrature approach is used to produce a numerical operator $N$ which approximates the switching operator of \eqref{eq_switch_op} as
\begin{equation}
 Nv^n(\br, q) = \frac{1}{\sqrt{\pi}} 
 \sum_{i=1}^M w_i v^n(r, \, \sqrt{2C}\sigma \xi_i + \theta + aC, \, 3-q). \label{eq_switch_hermite}
\end{equation}
To improve the accuracy, our implementation uses  $M = 3$ for this switching quadrature since $C \gg \tau$ in our experiments.

Our full method is summarized below in Algorithm \ref{alg_semilagrange}.
\begin{algorithm}
    \SetAlgoLined
    Initialize: set 
    $v^0 = 0$ on $\mathcal{T}$ and $v^0 = \infty$ elsewhere.\\
    Initialize: set $\texttt{maxdiff} = \infty$ and $n=0.$\\
    \While{$\texttt{maxdiff} > \epsilon$}{
        \ForEach{$r \in \mathcal{R}$}{
            \ForEach{$\theta \in \Theta$}{
                \ForEach{$q \in \{1, 2\}$}{
                    \ForEach{$u_i \in \mathcal{U}$} {
                        $\tau = 1.5 \min\left(\frac{\Delta r}{\dot{r}(q, u_i)}, \frac{\Delta \theta}{\dot{\theta}(r, q, u_i)}\right)$ \\
                        $S^\Delta_i = S_{u_i,\tau} v^n(1\br, q)$ \\
                    }
                $S^\Delta = \min_i S^\Delta_i$ \\
                $N^\Delta = C + N v^n(\br, q)$ \\
                $v^{n+1}(\br, q) = \min(S^\Delta, N^\Delta)$ \\
                }
            }
        }
        $\texttt{maxdiff} = \max_{\br, q}(|v^{n+1}(\br, q) - v^n(\br, q)|)$\\
        $ n = n+1$\\
    }
    \caption{ Semi-Lagrangian value iterations to a convergence criterion $\epsilon$. Here, $\mathcal{R}, \Theta, \mathcal{U}$ refer to the discrete set of gridpoints sampled for the state and control spaces.
}
    \label{alg_semilagrange}
\end{algorithm}

%\vspace{-0.5cm}

\subsection{Gauss-Seidel Iterations}
\label{ss:gauss_seidel}
% !TEX root = root.tex
Gauss-Seidel relaxation can greatly reduce the number of iterations needed to reach convergence. We heuristically find that time-parameterized stochastic-optimal paths $(r(t), \theta(t))$ almost always have monotone decreasing $r(t)$. Wherever this holds, the value function at any point will only depend on the values of gridpoints with smaller $r.$ This can be exploited by updating the value function estimates in-place, sweeping along increasing $r$ direction. As a result, grid points later in the $r$ sweep will see already-updated values at smaller $r$ rather than the old values from the previous iteration.

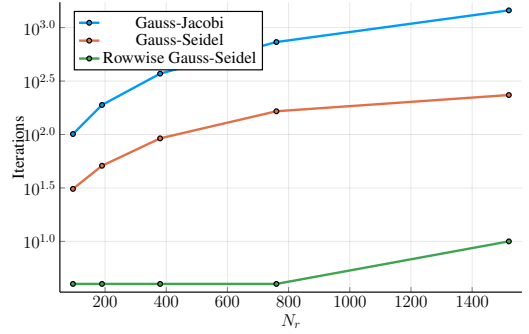
\begin{figure}[h!]
    \centering
    \resizebox{0.8\columnwidth}{!}{
        \begin{tikzpicture}[]
\begin{axis}[height = {101.6mm}, legend pos = {north west}, ylabel = {Iterations}, xmin = {52.25}, xmax = {1562.75}, ymax = {1730.4029330384674}, ymode = {log}, xlabel = {$N_r$}, unbounded coords=jump,scaled x ticks = false,xlabel style = {font = {\fontsize{14 pt}{18.2 pt}\selectfont}, color = {rgb,1:red,0.00000000;green,0.00000000;blue,0.00000000}, draw opacity = 1.0, rotate = 0.0},xmajorgrids = true,xtick = {200.0,400.0,600.0,800.0,1000.0,1200.0,1400.0},xticklabels = {$200$,$400$,$600$,$800$,$1000$,$1200$,$1400$},xtick align = inside,xticklabel style = {font = {\fontsize{14 pt}{18.2 pt}\selectfont}, color = {rgb,1:red,0.00000000;green,0.00000000;blue,0.00000000}, draw opacity = 1.0, rotate = 0.0},x grid style = {color = {rgb,1:red,0.00000000;green,0.00000000;blue,0.00000000},
draw opacity = 0.1,
line width = 0.5,
solid},axis x line* = left,x axis line style = {color = {rgb,1:red,0.00000000;green,0.00000000;blue,0.00000000},
draw opacity = 1.0,
line width = 1,
solid},scaled y ticks = false,ylabel style = {font = {\fontsize{14 pt}{18.2 pt}\selectfont}, color = {rgb,1:red,0.00000000;green,0.00000000;blue,0.00000000}, draw opacity = 1.0, rotate = 0.0},log basis y=10,ymajorgrids = true,ytick = {10.0,31.622776601683793,100.0,316.22776601683796,1000.0},yticklabels = {$10^{1.0}$,$10^{1.5}$,$10^{2.0}$,$10^{2.5}$,$10^{3.0}$},ytick align = inside,yticklabel style = {font = {\fontsize{14 pt}{18.2 pt}\selectfont}, color = {rgb,1:red,0.00000000;green,0.00000000;blue,0.00000000}, draw opacity = 1.0, rotate = 0.0},y grid style = {color = {rgb,1:red,0.00000000;green,0.00000000;blue,0.00000000},
draw opacity = 0.1,
line width = 0.5,
solid},axis y line* = left,y axis line style = {color = {rgb,1:red,0.00000000;green,0.00000000;blue,0.00000000},
draw opacity = 1.0,
line width = 1,
solid},    xshift = 0.0mm,
    yshift = 0.0mm,
    axis background/.style={fill={rgb,1:red,1.00000000;green,1.00000000;blue,1.00000000}}
,legend style = {color = {rgb,1:red,0.00000000;green,0.00000000;blue,0.00000000},
draw opacity = 1.0,
line width = 1,
solid,fill = {rgb,1:red,1.00000000;green,1.00000000;blue,1.00000000},fill opacity = 1.0,text opacity = 1.0,font = {\fontsize{14 pt}{18.2 pt}\selectfont}},colorbar style={title=}, ymin = {3.351820486004147}, width = {152.4mm}]\addplot+ [color = {rgb,1:red,0.00000000;green,0.60560316;blue,0.97868012},
draw opacity = 1.0,
line width = 2,
solid,mark = *,
mark size = 2.0,
mark options = {
            color = {rgb,1:red,0.00000000;green,0.00000000;blue,0.00000000}, draw opacity = 1.0,
            fill = {rgb,1:red,0.00000000;green,0.60560316;blue,0.97868012}, fill opacity = 1.0,
            line width = 1,
            rotate = 0,
            solid
        }]coordinates {
(95.0, 101.0)
(190.0, 189.0)
(380.0, 371.0)
(760.0, 733.0)
(1520.0, 1450.0)
};
\addlegendentry{Gauss-Jacobi}
\addplot+ [color = {rgb,1:red,0.88887350;green,0.43564919;blue,0.27812294},
draw opacity = 1.0,
line width = 2,
solid,mark = *,
mark size = 2.0,
mark options = {
            color = {rgb,1:red,0.00000000;green,0.00000000;blue,0.00000000}, draw opacity = 1.0,
            fill = {rgb,1:red,0.88887350;green,0.43564919;blue,0.27812294}, fill opacity = 1.0,
            line width = 1,
            rotate = 0,
            solid
        }]coordinates {
(95.0, 31.0)
(190.0, 51.0)
(380.0, 92.0)
(760.0, 165.0)
(1520.0, 234.0)
};
\addlegendentry{Gauss-Seidel}
\addplot+ [color = {rgb,1:red,0.24222430;green,0.64327509;blue,0.30444865},
draw opacity = 1.0,
line width = 2,
solid,mark = *,
mark size = 2.0,
mark options = {
            color = {rgb,1:red,0.00000000;green,0.00000000;blue,0.00000000}, draw opacity = 1.0,
            fill = {rgb,1:red,0.24222430;green,0.64327509;blue,0.30444865}, fill opacity = 1.0,
            line width = 1,
            rotate = 0,
            solid
        }]coordinates {
(95.0, 4.0)
(190.0, 4.0)
(380.0, 4.0)
(760.0, 4.0)
(1520.0, 10.0)
};
\addlegendentry{Rowwise Gauss-Seidel}
\end{axis}

\end{tikzpicture}
    }
    \caption{The number of iterations required to reach
    convergence with a fixed $\epsilon = 10^{-8}$ and an increasing number of grid points along each state dimension. 
    In this plot, we have fixed $\sigma = 0.05, a = 0.1$. At each point, $N_\theta$ is scaled accordingly with $N_r$.}
    \label{fig_sweep_perf}
\end{figure}

An additional change can be made to our discretization to further take advantage of this causal dependence. 
To ensure the currently updated point only depends on those at smaller $r$, we choose the adaptive timestep so that we always step at least a full grid length along the $r$ axis:

\begin{equation}
    \tau(\br, q, u) = 1.5 \frac{\Delta r}{\dot{r}(\br, q, u)}. \label{eq_row-wisedt}
\end{equation}
We refer to this second discretization as ``row-wise Gauss-Seidel" (rwGS) relaxation.
Our experiments show that it requires a nearly-constant number of iterations, only slowly beginning to increase at really fine grids. Meanwhile, both the Gauss-Jacobi (GJ) and the standard Gauss-Seidel (GS) require a number of iterations roughly linear in $N_R$ (the number of grid points along the $r$ direction).  As shown in Fig. \ref{fig_sweep_perf}, on the finest tested grid, GJ requires a factor of $\approx100$ more iterations than GS, which itself requires a factor of $\approx100$ more iterations than the rwGS.
However, this massive speedup in the latter is not completely free. Comparing \eqref{eq_adaptivedt} and \eqref{eq_row-wisedt}, we see that the row-wise iteration uses larger timesteps than the standard discretization in areas of the domain where $\Delta \theta / \dot{\theta} < \Delta r / \dot{r}$, in particular close to the target. Additional numerical tests (summarized in \href{https://eikonal-equation.github.io/Stochastic_Sailing/}{Supplementary Materials}) show that, despite slightly larger discretization errors in rwGS, these errors decrease with the same rate as in GS under the grid refinement.  

\subsection{Live Simulations / Control Synthesis}
\label{ss:sim_control}
% !TEX root = root.tex
Using the grid approximation of the value function obtained above, we can perform live sailboat navigation using two derived data structures. Our ``switchgrid'', indicates the gridpoints for which it was optimal to switch to the opposite tack (i.e., wherever $N^{\Delta} < S^{\Delta}$ in Algorithm \ref{alg_semilagrange}). We can also similarly keep track of a ``direction grid,'' recording the controls $u_i$ that minimize $S^{\Delta}$ at each gridpoint.

In our simulations, we evolve the system  \eqref{eq_reducedcoords} using a simple Euler–Maruyama method, with the stochastic component of the $\theta$ evolution sampled from $\mathcal{N}(0, \sigma^2)$ in each timestep. We switch to the opposite tack whenever at least $3$ of the surrounding $4$ gridpoints have $N^{\Delta} < S^{\Delta}$;
otherwise, we stay on the current tack and use the optimal control value of the nearest gridpoint in the direction grid.

\section{Numerical Experiments}
\label{s:experiments}

% !TEX root = root.tex
We now apply our algorithm to solve the optimal control problem \eqref{eq_reducedcoords}-\eqref{eq_valuefunc} for a variety of system parameters, and sample stochastic optimal trajectories for each scenario\footnote{
To ensure the computational reproducibility, we provide the full source code, movies illustrating these simulations, and additional convergence data at 
\url{https://eikonal-equation.github.io/Stochastic_Sailing/}
}. 
For all results shown, the target radius is $R = 0.1$ and the domain is $(r, \theta) \in [R, 2.0] \times [0, 2\pi]$. Target boundary conditions are set as $v(r, \theta, q) = 0, \; \forall r < R$ and exit boundary conditions are $v(r, \theta, q) = 10^6, \; \forall r > 2.0.$

\begin{figure*}[h!]
    \centering
    \begin{subfigure}{0.25\textwidth}
        \centering
        \includegraphics[width=\linewidth]{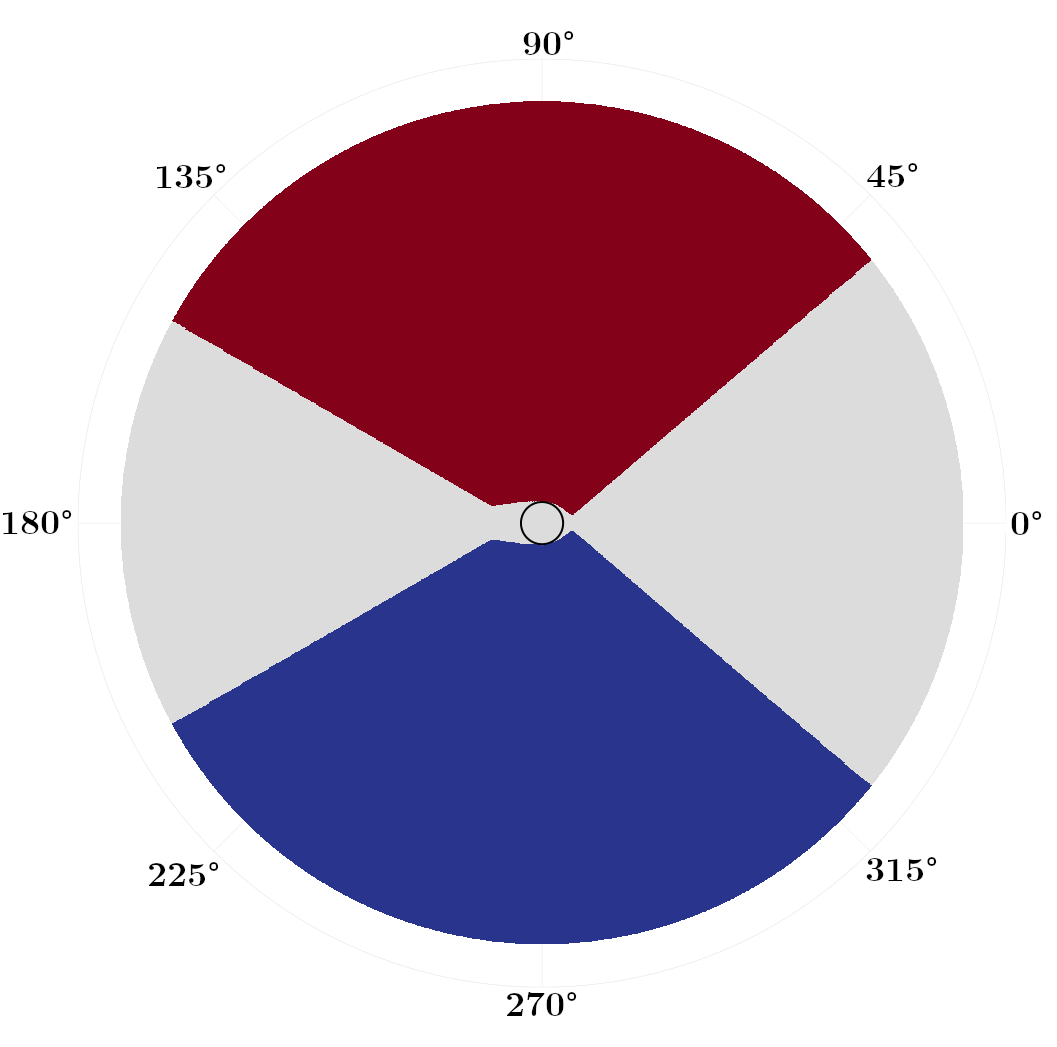}
        \caption{$\sigma = 0.0$}
    \end{subfigure}
    \hspace{1.1cm}
    \begin{subfigure}{0.25\textwidth}
        \centering
        \includegraphics[width=\linewidth]{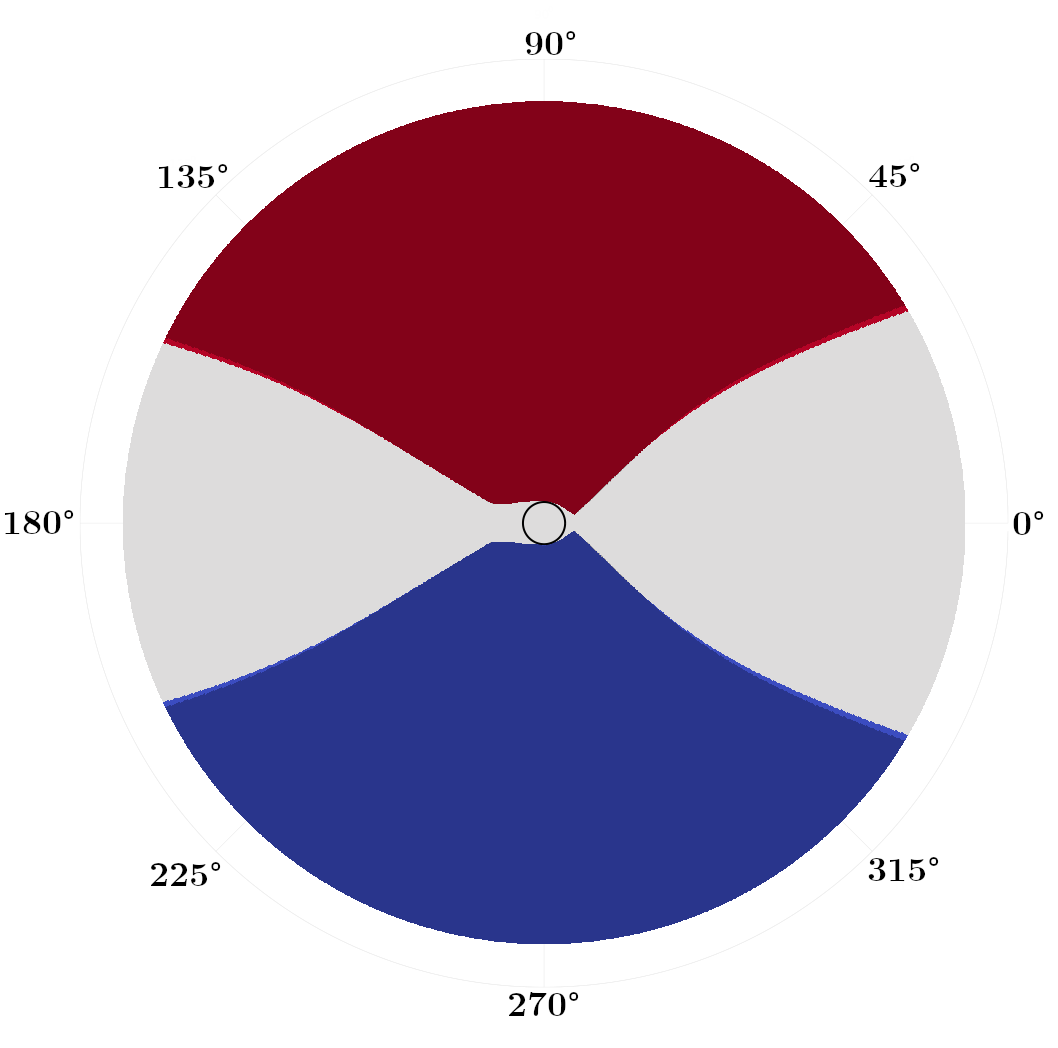}
        \caption{$\sigma = 0.05$}
    \end{subfigure}
    \hspace{1.1cm}
    \begin{subfigure}{0.25\textwidth}
        \centering
        \includegraphics[width=\linewidth]{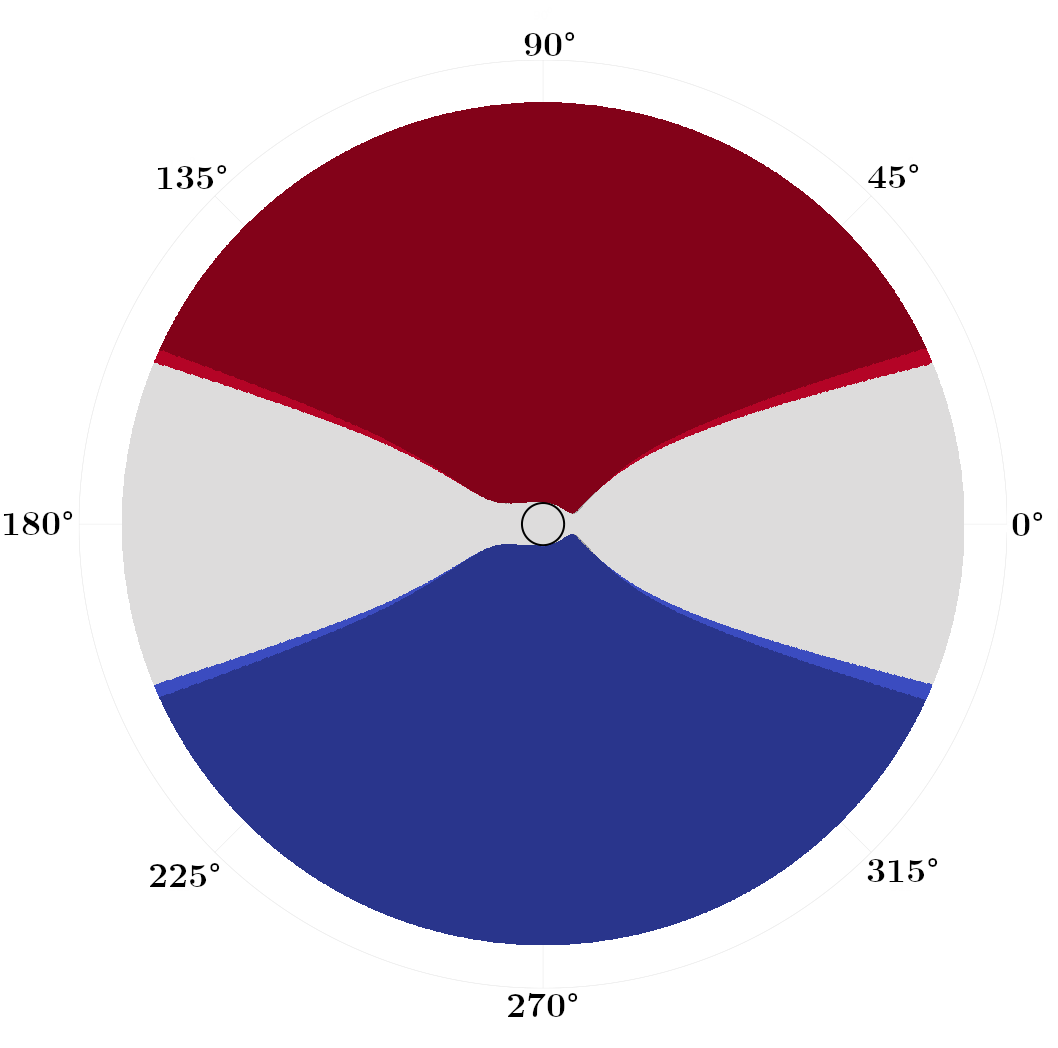}
        \caption{$\sigma = 0.10$}
    \end{subfigure}
    \caption{Switchgrids for stochastic zero-drift ($a = 0$) problems in $(r, \theta)$ space. Red points indicate that $\texttt{switchgrid}(r, \theta, q=1) = 1$, blue points indicate that $\texttt{switchgrid}(r, \theta, q=2) = 1$, and white points indicate that both are $0$. Darkened regions correspond to the switchgrids obtained when thezero-evolution switching operator $\mathcal{N} v(\br, q) = v(\br, 3 - q)$ from \cite{ferretti_optimal_2019} is used instead of our \eqref{eq_switch_hermite}.}
    \label{fig_sto_switch}
\end{figure*}

We first focus on problems with zero drift $(a=0)$, but varying
the diffusivity coefficient 
$\sigma$. In Fig.\ \ref{fig_sto_switch}, we plot switchgrid solutions for a collection of $\sigma$ values, noting that these grids are defined and plotted in $(r, \theta)$ state space, not the absolute position space.
Points are marked red if  $\texttt{switchgrid}(r, \theta, q=1) = 1$, blue if $\texttt{switchgrid}(r, \theta, q=2) = 1$, and white if both are $0$. The darker shaded regions correspond to the switchgrids obtained if we use the zero-evolution switching operator of \cite{ferretti_optimal_2019}
as compared to our \eqref{eq_switch_op}. The black circle at the center marks the target set $\mathcal{T}$. In general, stochastic optimal trajectories ``bounce'' between these red and blue switching fronts as they approach $\mathcal{T}$. 

\begin{figure*}[h!]
    \centering
    \begin{subfigure}{0.295\textwidth}
        \centering
        \includegraphics[width=\linewidth]{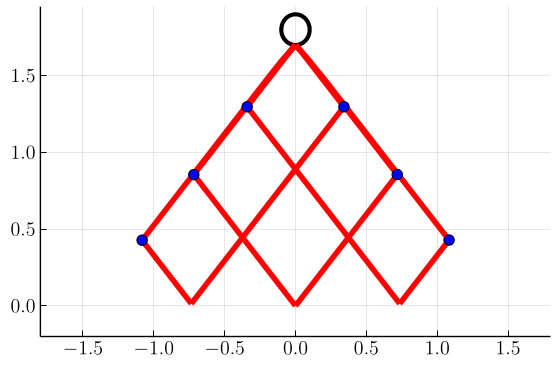}
        \caption{$\sigma = 0.0$}
    \end{subfigure}
    \begin{subfigure}{0.295\textwidth}
        \centering
        \includegraphics[width=\linewidth]{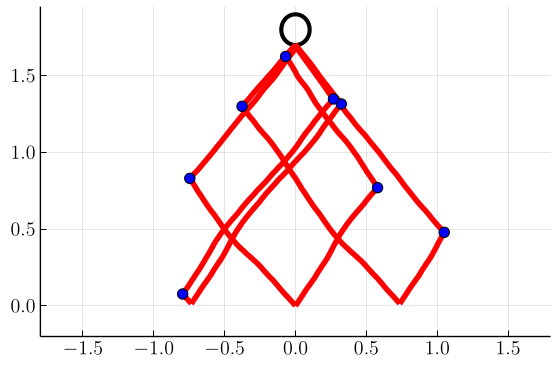}
        \caption{$\sigma = 0.05$}
    \end{subfigure}
    \begin{subfigure}{0.295\textwidth}
        \centering
        \includegraphics[width=\linewidth]{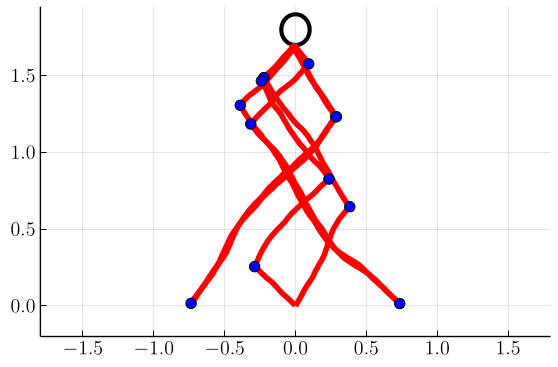}
        \caption{$\sigma = 0.10$}
    \end{subfigure}
    \caption{Six sample trajectories for  stochastic zero-drift ($a = 0$) problems: starting from 3 different $(r, \theta)$ positions and 2 different starting tacks ($q=1,2$). 
    Within each panel, all trajectories use the same sampled wind evolution. Blue points indicate the tack switch events.}
    \label{fig_sto_trajs}
\end{figure*}

With $\sigma = 0,$ the switching fronts lie at constant $\theta$, corresponding to the angles that locally maximize the projection of the polar speed plot (Fig.~\ref{fig_polarplot}) onto the axis of the wind. Once a sailboat hits this front, it can switch tack and stay at this constant angle that maximizes its speed towards the target (see Fig.~\ref{fig_sto_trajs}(a)).
As $\sigma$ increases,
these switching fronts contract more tightly around $\theta = 0^\circ, 180^\circ$. As the future wind state grows more uncertain with larger $\sigma$, the optimal policy becomes more conservative and constrains the boat to smaller excursions in $\theta$ before switching to the opposite tack. 
We observe that in these zero-drift problems, our improved switching operator makes the switching grids modestly more conservative, switching earlier by a few degrees in Fig.~\ref{fig_sto_switch}(b,c).

In Fig.~\ref{fig_sto_trajs}, we plot samples of stochastic optimal trajectories corresponding to these control policies, with blue markers indicating where tack switches occur. Similarly to \cite{ferretti_optimal_2019}, we initialize six boats at positions $(r, \theta) \in \{(1.8, 0.0), (1.93, \pm0.39)\},$ and tacks $q \in \{1, 2\}$. We then evolve and control each boat's state until it reaches its target as described in Sec. \ref{ss:sim_control}. We can see that, as foreshadowed by the switchgrids, at higher values of $\sigma$ the trajectories tend to stay within tighter cones from the target. 
\begin{figure*}[h!]
    \centering
    \begin{subfigure}{0.25\textwidth}
        \centering
        \includegraphics[width=\linewidth]{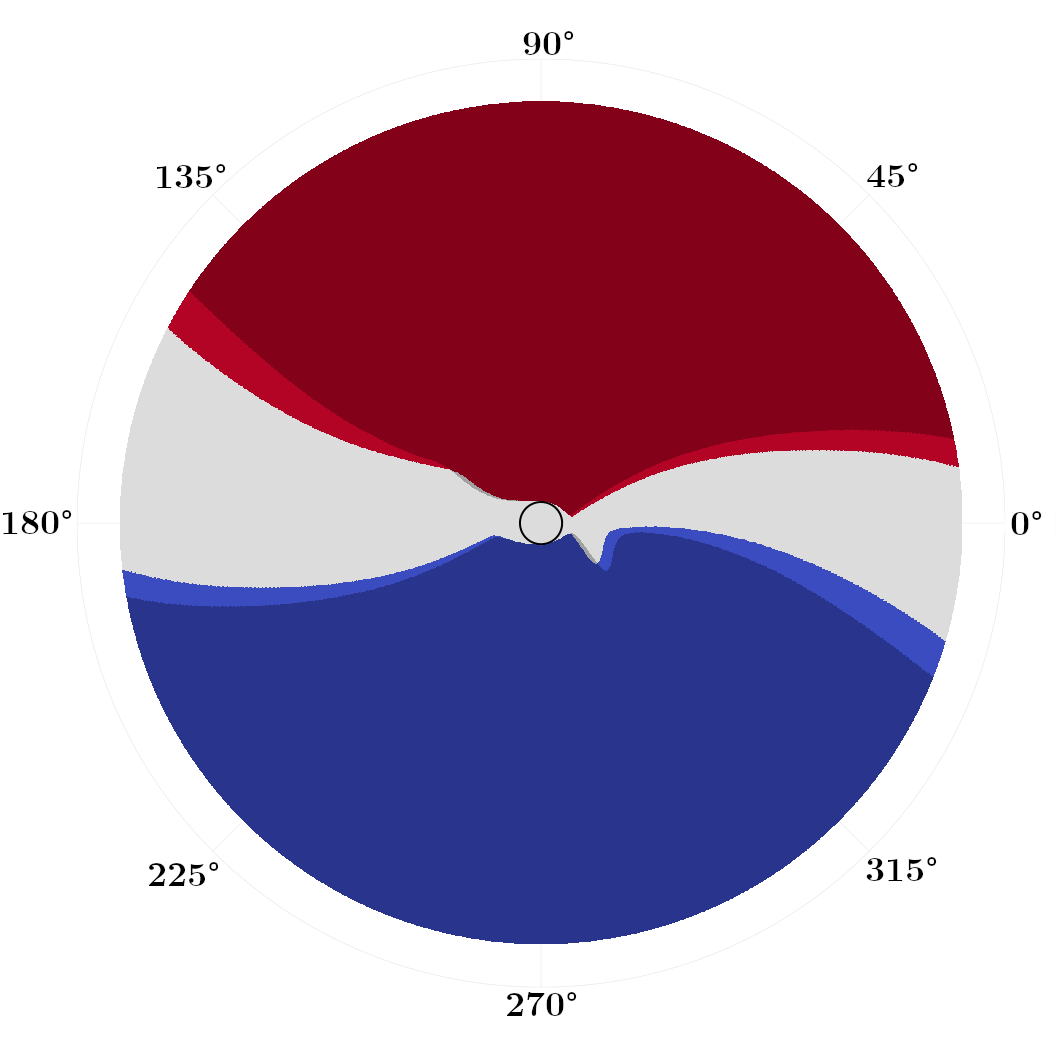}
        \caption{$a = 0.05$}
    \end{subfigure}
    \hspace{1.1cm}
    \begin{subfigure}{0.25\textwidth}
        \centering
        \includegraphics[width=\linewidth]{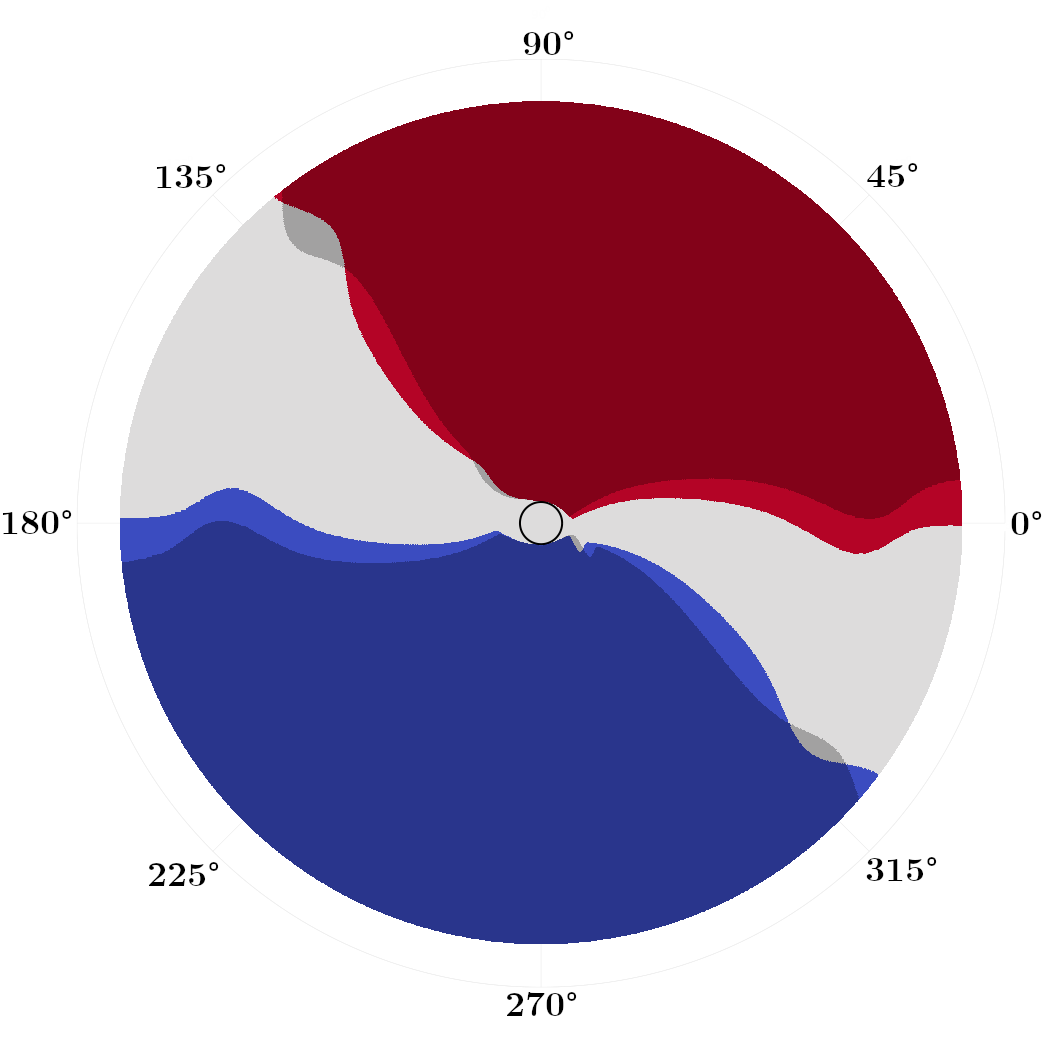}
        \caption{$a = 0.10$}
    \end{subfigure}
    \hspace{1.1cm}
    \begin{subfigure}{0.25\textwidth}
        \centering
        \includegraphics[width=\linewidth]{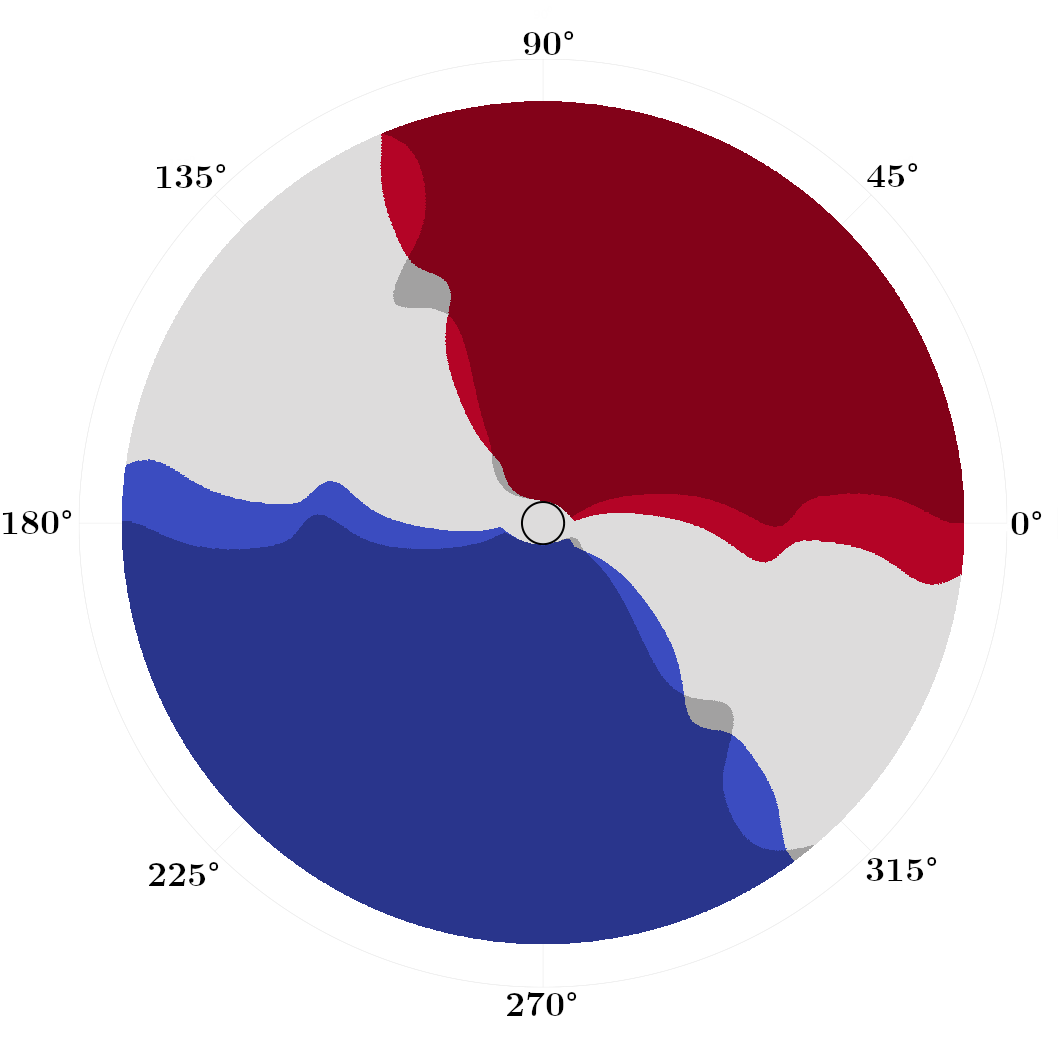}
        \caption{$a = 0.15$}
    \end{subfigure}
    \caption{Switchgrids for stochastic ($\sigma = 0.05$) nonzero-drift problems. Colorings are as in Fig.~\ref{fig_sto_switch}.}
    \label{fig_drift_switch}
\end{figure*}
\begin{figure*}[h!]
    \centering
    \begin{subfigure}{0.295\textwidth}
        \centering
        \includegraphics[width=\linewidth]{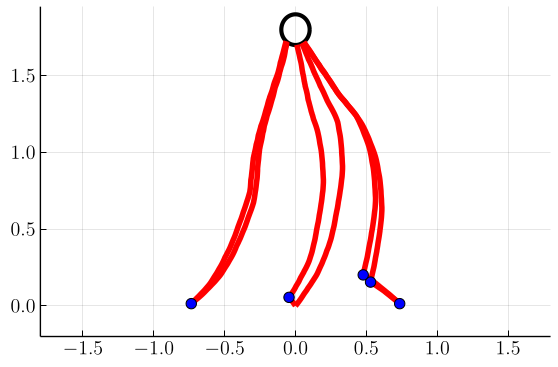}
        \caption{$a = 0.05$}
    \end{subfigure}
    \begin{subfigure}{0.295\textwidth}
        \centering
        \includegraphics[width=\linewidth]{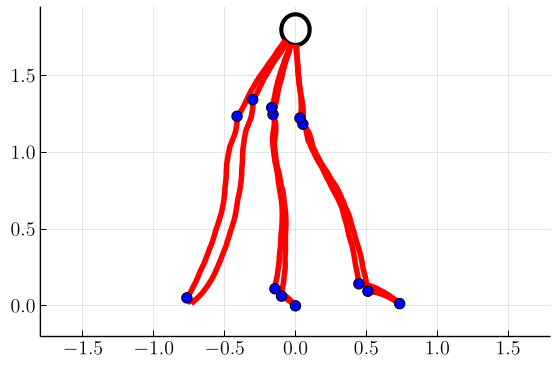}
        \caption{$a = 0.10$}
    \end{subfigure}
    \begin{subfigure}{0.295\textwidth}
        \centering
        \includegraphics[width=\linewidth]{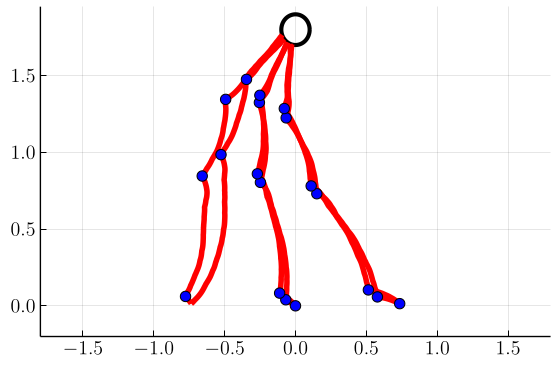}
        \caption{$a = 0.15$}
    \end{subfigure}
    \caption{Six sample trajectories for  stochastic ($\sigma = 0.05$) non-zero-drift problems: starting from 3 different $(r, \theta)$ positions and 2 different starting tacks.
    Within each panel, all trajectories use the same sampled wind evolution. Blue points indicate the tack switch events.}    
    \label{fig_drift_trajs}
\end{figure*}
Next, in Figs.\ \ref{fig_drift_switch} and \ref{fig_drift_trajs} we show switchgrids and sampled trajectories for problems in which the drift parameter $a$ is nonzero, with a fixed $\sigma=0.05$. As the drift increases, we observe the switchgrids and the trajectories becoming skewed and asymmetric towards negative $\theta$. This reflects that it is optimal to ``lead'' the wind in anticipation of future wind evolution, allowing the boat to stay in the portion of state space that allows it to steer at the optimal $~40^\circ \mhyphen 45^\circ$ angle while the wind evolves.

Additionally, some ``wave''-like features start to emerge in the switchgrids for large $a$. Empirically, we observe that the number of ``waves'' in the switchgrid corresponds to the number of $180^\circ$ turns the wind tends to make by the time an average boat at $r=R$ reaches the target. Each ``wave'' can thus be interpreted as corresponding to the optimal control during one (average) half-cycle of the wind. We notice that accurately modeling the switching operator is critical in these high-drift scenarios, as the two options produce dramatically different wave-like structures in the switchgrids.

\section{Conclusions}
\label{s:conclusions}
% !TEX root = root.tex

We have introduced an accelerated semi-Lagrangian approach to solve a hybrid stochastic optimal control problem arising in match race sailing. By working in a reduced coordinate system and designing a nearly-causal discretization, we were able to produce control policies on high-resolution grids which provide insight into tack-switching strategies under unpredictable weather.

While the conversion to reduced coordinates greatly accelerates our solution, it also introduces several limitations since we had to assume that the problem is radially symmetric. This restricts the target shape to a circle, and does not allow the inclusion of state constraints such as obstacles or nearby coastline.  
However, it does still allow for modeling match races between multiple agents.  Our approach will require keeping track of $(r,\theta)$ coordinates for each agent, yielding an overall state space with one less dimension than in \cite{cacace_stochastic_2019}.
For this reason, we hope that our approach will allow to rapidly test the implications of various weather models or match strategies in ``ideal'' conditions far away from any barriers. 

This paper is focused on a risk-neutral task of finding a feedback policy $\mu$ that minimizes the expected time to target $\mathbb{E}\left[T_\mu\right]$.  However, practitioners often prefer {\em robust} control policies more suitable when bad outcomes are costly.  Two such popular approaches are the ``risk-sensitive'' control \cite{HowardMatheson_1972, fleming1995risk} and $H^\infty$ control \cite{bacsar2008h}.  But neither of these provides any direct bounds on the likelihood of bad scenarios; i.e., $w_\mu(c) = \mathbb{P} \left( T_\mu \geq c \right).$
An efficient method for minimizing $w_\mu(c)$ for all threshold values $c$ simultaneously was recently introduced for piecewise-deterministic Markov processes in \cite{CarteeVlad_UQ}.
We hope that it can be similarly extended to the current context.  It would be also interesting to extend the methods of multiobjective optimal control (e.g., \cite{KumarVlad, DesillesZidani}) to find Pareto-optimal policies 
reflecting the rational tradeoffs between conflicting optimization criteria (e.g., $\mathbb{E}\left[T_\mu\right]$ vs $w_\mu(c)$ or even $w_\mu(c_1)$ vs $w_\mu(c_2)$).

\vspace*{3mm}
\noindent
{\bf Acknowledgements: } The authors are grateful to Roberto Ferretti for sparking their interest in this class of problems.

\bibliographystyle{IEEEtran}
\bibliography{IEEEabrv,mybibfile}

\end{document}